\documentclass[english,leqno]{amsart}

\usepackage{amsfonts,amssymb,amsmath}
\usepackage{setspace}
\usepackage{xcolor}
\usepackage{hyperref}

\usepackage[bb=dsserif]{mathalpha}

\newcommand{\elle} {\mathcal{L}}
\newcommand{\acca} {\mathcal{H}}

\newcommand{\erre} {\mathbb{R}}
\newcommand{\erren} {{\erre^{n}}}
\newcommand{\erreu} {{\erre^{{n + \scriptscriptstyle{1}}}}}

\newcommand{\frecciaf} {\longmapsto}
\newcommand{\meno} {\setminus}

\newcommand\intpalla{\int_{B(x_0,r)}}
\newcommand\intomegapalla{\int_{\Omega_r(z_0)}}
\newcommand \cappa{K}
\newcommand \gi{G}
\newcommand \palla{{B(x_0,r)}}
 \newcommand \omegapalla {{\Omega_r(z_0)}}
  \newcommand \omegapallarho {{\Omega_\rho(z_0)}}
 \newcommand   \unor {\frac{1}{r}}
 \newcommand\andd{ \quad\mbox{ and } \quad }
\newcommand\inn{\mbox{ in }}
\newcommand\ellestar{\elle^\ast}
\newcommand\menou{^{-1}}
\newcommand \supp {\mathrm{supp}\ }
\newcommand \zstar {z^\ast }

\newcommand \tstar {t^\ast }
\newcommand \ustar {u^\ast }
\newcommand \dstar {D^\ast }

\newtheorem{theorem}{Theorem}[section]

\newtheorem{corollary}[theorem]{Corollary}

\theoremstyle{remark}
\newtheorem{remark}[theorem]{Remark}

\theoremstyle{definition}

\def\stepone{\noindent{\it Step I. }}
\def\steptwo{\noindent{\it Step II. }}
\def\stepthree{\noindent{\it Step III. }}

\numberwithin{equation}{section}

\title[A rigidity theorem for Kolmogorov-type operators]
 {A rigidity theorem \\for Kolmogorov-type operators}
 
 \author{Alessia E. Kogoj}
\address{Dipartimento di Scienze Pure e Applicate (DiSPeA)\\ 
				 Universit\`{a} degli Studi di Urbino Carlo Bo\\
				 Piazza della Repubblica 13, 61029 Urbino (PU), Italy.}
\email{alessia.kogoj@uniurb.it}

\author{Ermanno Lanconelli}
\address{Dipartimento di Matematica\\ 
				 Alma Mater Studiorum Universit\`{a} di Bologna\\
				 Piazza di Porta San Donato 5, 40126 Bologna, Italy.}
\email{ermanno.lanconelli@unibo.it}



\begin{document}
\maketitle
{\small \emph{To Enzo Mitidieri, on the occasion of
his 70th birthday, \\ with profound friendship and great admiration.}}

\begin{abstract}
Let $D\subseteq \erren$, $n\geq 3$, be a bounded open set and let $x_0\in D$. Assume that the Newtonian potential of $D$ is 
proportional outside $D$ to the Newtonian potential of a mass concentrated at $\{x_0\}.$  Then $D$ is a Euclidean ball centered at $x_0$. This Theorem, proved by Aharonov, Shiffer and Zalcman in 1981, was extended to the caloric setting by Suzuki and Watson in 2001. In this note, we show that Suzuki--Watson Theorem is a particular case of a more general rigidity result related to a class of Kolmogorov-type PDEs.

\end{abstract}

\maketitle

\section{Introduction}
\subsection{Rigidity results in harmonic and caloric settings}\mbox{}

Let $\cappa$ be the fundamental solution with pole at the origin of the Laplacian $\varDelta$ in $\erren$, $n\geq 3$.
We denote by $\palla$ the Euclidean ball of $\erren$ centered at $x_0$ with radius $r>0$. For every $y \notin \palla$ 
the function $x\frecciaf \cappa(x-y)$ is harmonic in $\palla$ so that, by Gauss Mean Value Theorem, 

\begin{equation}\label{uno} \intpalla \cappa(x-y)\ dx = c\, \cappa(x_0-y),\end{equation}
where $c$ is the volume of $\palla$.

In 1981, Aharonov, Shiffer and Zalcman proved that identity \eqref{uno} is a rigidity property of the Euclidean ball. More precisely, they proved (see \cite{potato}, see also \cite{cuplanidea}) that a bounded open set $D$ such that, for a point $x_0\in D$ and a suitable positive constant $c$, 

\begin{equation*}\ \int_D \cappa(x-y)\ dx = c\, \cappa(x_0-y) \quad \forall\, y\notin D,\end{equation*}
has to be the Euclidean ball $B(x_0,r)$ with Lebesgue measure equal to $c$.

Suzuki and Watson, in 2001, extended the previous theorem to the {\it heat balls} in $\erreu$, $n\geq 1$  (see \cite{suzuki_watson}). To be more precise, we need some notation. 

Let us denote by $\gi$ the fundamental solution of the heat operator $\acca:=\varDelta -\partial_t$ in 
$\mathbb R^{n+1} = \mathbb R^n_x \times \mathbb R_t$. We call {\it heat ball}  with center $z_0\in \erreu$ and $r>0$ the following set 
$$
 \Omega_r(z_0):  = \left\{z \in \mathbb R^{n+1} : \gi (z_0 - z) > \frac{1}{r} \right\}.
$$

The caloric functions, i.e., the solutions to the heat equation 
$$\acca u=0$$ 
can be characterized in terms of the following caloric mean value formula 

\begin{equation}\label{cmvf} u(z_0)=M_r(u)(z_0):= \frac{1}{r} \intomegapalla u(\zeta) W(\zeta-z_0)\ d\zeta, \end{equation}
where $W$ is the Pini--Watson kernel 
\begin{equation}\label{ast} W(\eta)=W(\xi,\tau):= \frac{1}{4} \left( \frac{|\xi|}{\tau}\right)^2, (\xi,\tau)\in\erreu,\tau\neq 0. \end{equation}

Indeed, a continuous function $u:O\longrightarrow \erre,$ $O$ open subset of $\erreu$, is smooth and solves $\acca u=0$ in $O$ if and only if 

$$u(z_0)= M_r(u)(z_0)$$
for every $z_0\in O$ and $r>0$ such that $\overline{\omegapalla} \subseteq O$ (see \cite{W}).

As a consequence, for every heat ball $\omegapalla$ and for every $z\notin\omegapalla$, one has 

\begin{equation}\label{due} \frac{1}{r} \int_{\omegapalla} G(\zeta-z) W(\zeta-z_0)\ d\zeta =G(z_0-z).\end{equation}

Indeed, if $0<\rho<r$ and and $z\notin \omegapalla$, $z\neq z_0$, then 
$$\zeta \frecciaf G(\zeta- z)$$ is caloric in $\erreu\meno \overline\omegapallarho.$ Hence, by the caloric mean value formula \eqref{cmvf},

\begin{equation*} \frac{1}{\rho}\int_\omegapallarho G(\zeta -z) W(\zeta-z_0)\ d\zeta = G(z_o-z).\end{equation*}


From this identity, as $\rho\nearrow r$, one gets \eqref{due} in the case $z\neq z_0.$ On the other hand, if $z=z_0$ 
identity  \eqref{due}  is trivial. Suzuki and Watson, extending Aharonov, Shiffer and Zalcman's Theorem to the caloric setting, proved that \eqref{due}  is a rigidity property of the heat balls.

Their result reads as follows: let $z_0=(x_0,t_0)\in\erreu$ and let $D$ be a bounded open set of $\erreu$. Assume that for a suitable constant $c>0$, 

\begin{equation*} \int_D G(\zeta -z) W(\zeta-z_0)\ d\zeta = c\, G(z_o-z)\quad \forall\, z\notin D.\end{equation*}

Then, if

\begin{equation}\label{tre} \zeta \frecciaf (\mathbb{1}_D -\mathbb{1}_{\Omega_c(z_0)})(\zeta) W(z_0-\zeta)\in L^p \mbox{\ for some\ } p>\frac{n}{2}+1,\end{equation}
then $$D=\Omega_c(z_0).$$

Here, and in what follows,  $\mathbb{1}_E$ denotes the characteristic function:
  $\mathbb{1}_E(x)=1$ if $x\in E$,  $\mathbb{1}_E(x)=0$ otherwise.

We want to stress that condition  \eqref{tre} replaces the condition $x_0\in D$ of the harmonic case; its meaning is that $D$ and $\Omega_c(z_0)$ are indistinguishable in the vicinity of $z_0$.
The present authors, together with G. Tralli, in \cite{inverse} partially extended this last rigidity theorem to a class of second order hypoelliptic operators containing in particular the prototype of the so called Kolmogorov operators.

 In this note, we provide a full extension of Suzuki--Watson's rigidity Theorem to such a class of partial differential operators. Our technique is inspired by the one used in the paper \cite{cuplanidea} where 
harmonic characterizations of the Euclidean balls are proved. 
\subsection{Our operators}\mbox{}

We will deal with Partial Differential Operators of  the type:
\begin{eqnarray}\label{KFP}
\mathcal L 
: = \mathrm{div}\left(  A\nabla\right)  +\left\langle Bx,\nabla\right\rangle - \partial_t,
\end{eqnarray} where $(x,t) \in \mathbb R^n \times \mathbb R$, $\nabla$ and $\left\langle \cdot,\cdot \right\rangle$ denote the gradient and the inner product in $\erren$. $A = (a_{i,j})_{i,j = 1, \dots , n}$ and $B= (b_{i,j})_{i,j = 1, \dots , n}$ are $n\times n$ real constant matrices taking the following block form: 

\begin{equation*}
A=%
\begin{bmatrix}
A_{0} & 0\\
0 & 0
\end{bmatrix} \label{A}%
\end{equation*}\\
\begin{equation*}
B=%
\begin{bmatrix}
0 & 0 &  \ldots & 0 & 0 \\
B_1& 0 &  \ldots & 0  & 0 \\
0& B_2&  \ldots & 0  & 0 \\
\vdots & \vdots & \ddots & \vdots & \vdots\\
0 & 0 &  \ldots & B_{r}&0
\end{bmatrix},
\label{B}%
\end{equation*}\\
where $A_0$ is a  $p_{0}\times p_{0}$ ($1\le p_{0}\leq n$) symmetric and positive definite
constant matrix and $B_{j}$ is a $p_{j}\times p_{j-1}$ ($j=1,2,...,r$)  block with rank equal to $p_{j}$. Moreover $p_{0}\geq p_{1}\geq...\geq p_{r}\geq1$ and $p_{0}+p_{1}+...+p_{r}=n$.

We explicitly remark that the operator \eqref{KFP} becomes the heat operator if $A=\mathbb{I}_n$ - the identity matrix - and 
$B=0$. In this case, with the previous notations, $p_0=n$ and $p_1,\ldots, p_r$ disappear.

It is quite well known that, under these block form assumptions on $A$ and $B$, the operator $\elle$ in \eqref{KFP}  is hypoelliptic, i.e., every distributional solution $u$ to $\elle u=f$ is actually smooth whenever $f$ is smooth 
(see \cite{lanconelli_polidoro_1994},  see also \cite[Chapter 4, Section 4.3.4]{BLU}).
It is also well known that $\elle$  is left translation invariant and homogeneous of degree two on the homogeneous Lie group $$\mathbb {K}:=(\erreu, \circ, \delta_\lambda)$$ whose composition law is the following one
$$
\left(  x,t\right)  \circ\left(  \xi,\tau\right)  =\left(  \xi+E(\tau)  x,t+\tau\right),
$$ 
with $E(\tau) =\exp\left(-\tau B\right);$ moreover the dilation $\delta_\lambda, \lambda >0$,  is  the linear map from 
$\erreu$ to $\erreu$ whose Jacobian matrix is given by 
$$D(\lambda):=\mathrm{diag} (\lambda \mathbb{I}_{p_0}, \lambda^3  \mathbb{I}_{p_1}, \ldots, \lambda^{2r+1} \mathbb{I}_{p_r}, \lambda^2),$$ being $\mathbb{I}_{p_j}$ the $p_j\times p_j$ identity matrix. We remark that 
$$\det D(\lambda)=\lambda^Q,$$ with $Q:=p_0+p_1+\ldots + (2k+1)p_r+2.$
This natural number is the {\it homogeneous dimension of $\mathbb{K}$}. Since  $\mathbb{K}$ is a homogeneous Lie group in $\erreu$, the Lebesgue measure in $\erreu$ is left and right translation invariant on $\mathbb{K}$.

An explicit fundamental solution for  \eqref{KFP} is given by 

\begin{equation}\label{fundsol1}
\Gamma\left(  z,\zeta\right):  =\gamma\left(  \zeta^{-1}\circ z\right)  \text{
for }z,\zeta\in\mathbb{R}^{n+1},
\end{equation}
where $\zeta^{-1}=\left(  \xi,\tau\right)  ^{-1}=\left(  -E\left(  -\tau\right)  \xi,-\tau\right)$ denotes the opposite of $\zeta$ with respect to the composition law in $\cappa$ and 
\begin{equation}\label{fundsol2}
\gamma(z)=\gamma\left(x,t\right): =\left\{
\begin{tabular}
[c]{ll}%
$0$ & $\text{for }t\leq0$\\
$\frac{\left(  4\pi\right)  ^{-N/2}}{\sqrt{\det C\left(  t\right)  }}%
\exp\left(  -\frac{1}{4}\left\langle C^{-1}\left(  t\right)  x,x\right\rangle
\right)  $ & $\text{for }t>0$%
\end{tabular}
\ \right.. \end{equation}
Here $C\left(  t\right)  $ stands for the matrix 

\begin{equation*}
C(t):=\int_{0}^{t} E(s)AE^T(s)\,ds.  %
\end{equation*}

This matrix is strictly positive definite for every $t>0$ and strictly negative definite for every $t<0$ (see \cite{lanconelli_polidoro_1994}). In the case of the heat operator, $C(t)$ is simply given by $t\mathbb{I}_n.$

The function $\gamma$, the fundamental solution of $\elle$ with pole at the origin - the neutral element of $\mathbb{K}$ - is $\delta_\lambda$-homogeneous of degree $2-Q$, i.e., 

$$\gamma(\delta_\lambda (z)) = \lambda^{2-Q} \gamma(z)\quad \forall\, z\in\erreu,\ \forall\, \lambda>0.$$

\subsection{Mean Value formula}\mbox{}\\
Let $\elle$ be the operator \eqref{KFP}. For every $z_0 \in \erreu$ and $r>0$ we call {\it $\elle-$ball with center $z_0$ and radius $r>0$} the following open set
\begin{equation}\label{palle}
 \Omega_r(z_0) = \left\{z \in \mathbb R^{n+1} : \Gamma (z_0, z) > \frac{1}{r} \right\}.
\end{equation}

From \eqref{fundsol1} and \eqref{fundsol2}, one easily verifies that $\omegapalla$ is a non-empty bounded open set; moreover,

$$\bigcap_{r > 0} {\Omega_r(z_0)} = \{z_0\}.$$

A continuous function $u:O\longrightarrow \erre$, $O\subseteq \erreu$ open, actually is smooth in $O$ and solves 
$$\elle u=0 \inn O$$
if and only if 
\begin{equation}\label{sei} u(z_0)= \unor \intomegapalla u(\zeta) W(z_0^{-1} \circ \zeta)\ d \zeta,
\end{equation}
for every $\elle$-ball $\omegapalla$ such that $\overline\omegapalla\subseteq O$ (see \cite{cup_lan_media}, see also \cite[Theorem 1.1]{inverse}).  In \eqref{sei} the kernel $W$ is defined as follows

\begin{equation}\label{sette} W(z)=W(x,t):=\frac{\langle A C^{-1}(t) x, C^{-1}(t) x\rangle}{4}.
\end{equation}
$W$ is a well-defined and strictly positive almost everywhere in $\erreu$ smooth function. Indeed, since $A\geq 0$,  $W\geq 0$ in $\erreu \meno(\erren \times \{0\}).$ Moreover, as 
$$ C^{-1}(t)>0\quad\forall\, t>0 \andd C^{-1}(t)<0\quad\forall\, t<0,$$ 
one has $W(x,t)=0$ if and only if $x\in F_t:=C(t) (\ker (A)).$ 
Being $\mathrm{rank}(A)=p_0>0,$ $C(t) (\ker (A))$ has dimension $n-p_0$, hence strictly less then $n$. It follows that $F_t$
has $n$-measure equal to zero for every $t\neq 0$  and this implies that 
$$F:= \{ (x,t)\in \erreu \ |\ W(x,t)=0\}$$ has $n+1$-measure equal to zero.

We remark that when $\elle=\acca$, the kernel  \eqref{sette} becomes the Pini--Watson kernel in \eqref{ast}. 

From the Mean Value formula \eqref{sei}, just proceeding as in the caloric case, one gets 
\begin{equation}\label{otto}  \intomegapalla \Gamma (\zeta,z) W(z_0^{-1} \circ \zeta)\ d \zeta = r \Gamma(z_0,z) \quad \forall\, z\notin \omegapalla.
\end{equation}

\subsection{Main Theorem}\mbox{}\\
The aim of this paper is to prove that identity \eqref{otto} is a rigidity property of the $\elle$-balls; equivalently, we want to extend the Suzuki--Watson Theorem to the $\elle$-setting. Here is our main theorem in which $\Gamma$ denotes the fundamental solution of the operator $\elle$ in \eqref{KFP}, $\omegapalla$ is the $\elle$-ball in \eqref{palle} and $Q$ is the homogeneous dimension of $\mathbb{K}.$
\begin{theorem}\label{main} Let $z_0\in\erreu$ and let $D$ be a bounded open subset of $\erreu$ such that, for a suitable $r>0$, 
\begin{equation}\label{trestelle}  \int_D \Gamma (\zeta,z) W(z_0^{-1} \circ \zeta)\ d \zeta = r \Gamma(z_0,z) \quad \forall\, z\in \erreu\meno D.
\end{equation}
If, moreover, 
\begin{equation}\label{nove} (\mathbb{1}_D -\mathbb{1}_{\Omega_r(z_0)} )W(z_0^{-1} \circ \cdot)\in L^p \mbox{ for some } p>\dfrac{Q}{2},\end{equation}
then $D=\omegapalla.$

\end{theorem}

We explicitly remark that identity \eqref{trestelle} implies the inclusion $$D\subseteq \erren \times ]-\infty, t_0[,$$
being $t_0\in\erre$ the time-component of $z_0$, i.e., $z_0=(x_0,t_0)$ for a suitable $x_0\in\erren.$ Indeed, the right hand side of \eqref{trestelle} is equal to zero for every $z=(x,t)$ with $t\geq t_0$ while $\Gamma(\zeta, z)>0$ if $\zeta=(\xi,\tau)\in D$ and $\tau>t$.

Due to the invariance of the Lebesgue measure with respect to the right (and left) translation on $\mathbb{K}$, identity \eqref{trestelle}  is equivalent to the following one 
\begin{equation}\label{trestelle'}  \int_{z_0^{-1} \circ D}  \Gamma (\zeta,z) W(\zeta)\ d \zeta = r \Gamma(0,z) \quad \forall\, z\notin z_0^{-1}\circ D.
\end{equation}
Then, it is enough to prove Theorem \ref{main} in the case $z_0=0$.
\begin{corollary}\label{cor} Let $z_0\in\erreu$ and let $D$ be a bounded open subset of $\erreu$ such that, for a suitable $r>0$, 
\begin{equation}\label{star!}  u(z_0)= \unor \int_{D}  u(\zeta) W(z_0^{-1} \circ \zeta)\ d \zeta \end{equation}
for every non negative function $u$ $\elle$-harmonic in an open set containing $D\cup \{z_0\}.$
Then 
$$D=\omegapalla$$
if condition \eqref{nove} is satisfied. 
\end{corollary}
\begin{remark} If we replace condition \eqref{nove} with the following stronger ones:
\begin{itemize} \item[(i)] there exists a neighborhood $V$ of $z_0$ s. t. $\omegapalla \cap  V=D \cap V$;
         \item[(ii)]  $\overline{D}\meno \{z_0\}\subset \erren \times ]-\infty,t_0[,$  being $t_0\in\erre$ the time-component of $z_0$, i.e., $z_0 = (x_0, t_0)$ for a suitable $x_0\in \erren;$\end{itemize} 
then Theorem \ref{main} becomes a particular case of Theorem 1.4 in \cite{inverse}.
\end{remark}

The paper is organized as follows. In Section 2 we recall some notion and several results from Potential Analysis for the operator $\elle$ and its adjoint $\elle^\ast$. These will be the main ingredients and tools of our proof of Theorem  \ref{main}. Section 3 
will be devoted entirely to the proof of Theorem  \ref{main} and Section 4 to the proof of its Corollary \ref{cor}.

\section{Basic notions and results from Potential Analysis for $\elle$ and $\elle^\ast$}

\subsection{Harmonic and subharmonic functions for $\elle$ and $\elle^\ast$}\ \mbox{}\\
Let $\Omega\subseteq \erreu$ be open. A function $u:\Omega\rightarrow \erre$ is called {\it $\elle^\ast$-harmonic} -in short notation $u\in \elle^\ast(\Omega)$ - if $u\in C^\infty(\Omega,\erre)$ and $\ellestar u=0$ in $\Omega$. We explicitly observe that 
\begin{eqnarray*}
\ellestar
: = \mathrm{div}\left(  A\nabla\right)  -\left\langle Bx,\nabla\right\rangle + \partial_t,
\end{eqnarray*}
and therefore, as $\elle$, it is hypoelliptic.

\noindent Analogous meanings as above for 
$\elle$-harmonic functions and for $\elle(\Omega).$ 

A bounded open set $V\subseteq\erreu$ is called {\it $\ellestar$-regular} if for every function $\varphi\in C(\partial V,\erre)$ there exists a unique $\ellestar$-harmonic function in $V$, denoted by $H_\varphi^V$, such that 

$$\lim_{z\rightarrow\zeta}H_\varphi^V(z)=\varphi(\zeta)\quad\forall\, \zeta\in\partial V.$$
Analogous meaning for $\elle$-regular set.

Let $\Omega\subseteq\erreu$ be open and let $u:\Omega\longrightarrow [-\infty,+\infty[$ be an upper semicontinuous function.
We say that $u$ is $\ellestar$-subharmonic in $\Omega$ - in short notation $u\in\underline\elle^\ast (\Omega)$ - if it satisfies 
the following conditions:

\begin{itemize}
\item[(i)] $u>-\infty$ in a dense subset of $\Omega$;
\item[(ii)] for every $\ellestar$-regular open set $V$ with $\overline V\subseteq \Omega$ and for every $\varphi\in C(\partial V,\erre)$ such that $u|_{\partial V}\le \varphi$ one has $u\le H_\varphi^V$ in $V$.
\end{itemize}
We shall denote by $\overline\elle^\ast(\Omega)$ the family of the {$\ellestar$-superharmonic functions, i.e.,  the family of the functions $v$ such that $-v\in \underline\elle^\ast(\Omega).$
\subsection{Maximum principle for $\ellestar$-subharmonic functions}\label{subsection 2.3}\ \mbox{}\\
Let $\Omega\subseteq\erreu$ be a bounded open set and let $u\in\underline\elle^\ast(\Omega)$ be such that 
$$\limsup_{z\rightarrow \zeta} u(z)\le 0\quad\forall\, \zeta \in\partial\Omega.$$
Then
$$u\le 0\inn\Omega$$
(see e.g. \cite[Proposition 2.3]{cintilanconelliPJ2009}).

\subsection{Propagation of maxima along drift-trajectories}\label{subsection 2.4}\ \mbox{}\\
We call {\it drift-trajectory} of $\ellestar$ any path of the type
$$ s\frecciaf \gamma(s):=\alpha+s e_j \quad \mbox{ or } \quad s\frecciaf \gamma(s):=\alpha-s e_j,$$
where $0\le s\le S,\quad  \alpha \in \erreu, \quad e_j=(0,\ldots,0,\underset{j}{1},0,\ldots, 0),\quad 0\le j\le p_0.$ 

Then, if $u\in\overline\elle^\ast (\Omega), \Omega$ open subset of $\erreu$, and $z_0\in\Omega$ is such that 
$$u(z_0)=\max_\Omega u,$$
then $u(\gamma(s))=u(z_0)$ for every drift-trajectory $\gamma:[0,S]\longrightarrow\Omega$ such that $\gamma(0)=z_0$ 
(see e.g. \cite{CNP2010}, \cite{kogoj_polidoro}).

\subsection{$\Gamma$-potentials}\label{subsection 2.5} Let $\Gamma$ be the fundamental solution of $\elle$ defined in \eqref{fundsol1}. 
and $\mu$ be a non-negative Radon measure with compact support. We let $\Gamma_\mu: \erreu\longrightarrow [0,\infty],$
$$\Gamma_\mu(z):=\int_\erreu \Gamma(z,\zeta)\ d\mu(\zeta).$$
Then, see e.g. \cite[Proposition 4.1]{cintilanconelliPJ2009}, 
$$\Gamma_\mu\in\overline\elle^\ast(\erreu)$$
and 
$$\elle^\ast \Gamma_\mu=-\mu\quad \mbox{in the weak sense of the distributions.}$$
In particular $\Gamma_\mu$ is $\ellestar$-harmonic in $\erreu\meno \mathrm{supp}\mu$. We want to explicitly remark that  $\erreu\meno \mathrm{supp}\mu$ is the open set union of the family of the open sets $O$ such that $\mu(O)=0$.

\subsection{An inequality for the $\Gamma$-potentials of the $\elle$-balls}  \mbox{}\\
For fixed $z_0\in\erreu$ and $r>0$, denote by $\mu$ the Radon measure in $\erreu$ such that 
\begin{equation}\label{dieci} d\mu(\zeta)= \unor\mathbb{1}_\omegapalla (\zeta) W(z_0\menou\circ\zeta)\ d\zeta,
\end{equation}where  $\omegapalla$ is the $\elle$-ball centered at $z_0$ with radius $r$ and $W$ is the kernel \eqref{sette}. We have already noticed - see \eqref{otto} - that 
$$\Gamma_\mu(z)=\Gamma (z_0,z)\quad \forall\, z\notin \omegapalla.$$

From  Corollary 3.2 in \cite{kt}, we also get the following inequality 
\begin{equation}\label{undici} \Gamma_\mu(z)<\Gamma(z_0,z)\quad \forall\, z\in\omegapalla.
\end{equation}

This inequality will play a crucial r\^ole in the proof of the Theorem \ref{main}. Here, we stress that, together with \eqref{sette}, it implies

\begin{equation}\label{stellina} \Gamma_\mu(z)\le \Gamma(z_0,z)\quad \forall\, z\in\erreu.
\end{equation}

\subsection{A convolution continuity result}\label{subsection 2.7}\mbox{}\\
Since $\gamma$ is $\delta_\lambda$-homogenous of degree $2-Q$ we have 
$$\gamma\in L^q_{\mathrm{loc}}(\erreu) \mbox{\quad if } q(Q-2)<Q$$
or, equivalently, if 
$$0<\frac{1}{p}:=1 - \frac{1}{q}<\frac{2}{Q}.$$
Then, if $f\in L^p(\erreu)$ for some $p>\dfrac{Q}{2}$ and, moreover, the support of $f$ is compact, 
the function 
$$z\frecciaf \Gamma(z):=\int_{\erreu} \Gamma(\zeta,z) f(\zeta)\ d\zeta$$
is a well-defined real continuous function in $\erren$. The proof of this statement is completely standard if one remarks that, being the Lebesgue measure translation invariant in $\mathbb{K}$, 
$$\int_\erreu \Gamma(\zeta,z) f(\zeta)\ d\zeta = \int_\erreu \gamma(z\menou \circ \zeta)f(\zeta)\ d\zeta=\int_\erreu \gamma(\eta) f(z\circ \eta)\ d\eta.$$

\section{Proof of Theorem \ref{main}}
As we already noticed, it is not restrictive to assume $z_0=0$. Let $\mu$ be the compactly supported Radon measure defined in \eqref{dieci}, and let $\nu$ be the measure such that 
\begin{equation}\label{dodici} d\nu(\zeta)= \unor\mathbb{1}_D (\zeta) W(z_0\menou\circ\zeta)\ d\zeta.
\end{equation}
Since $D$ is bounded, $\nu$ is also a compactly supported Radon measure. For our aims, it is crucial to remark that 
$\mathrm{supp}\nu =\overline D$, hence 
\begin{equation}\label{tredici} \partial D \subseteq  \mathrm{supp} \nu.
\end{equation}
It is also crucial for us to remark that 
\begin{equation}\label{14}  \mu|_{\omegapalla\cap D} = \nu|_{\omegapalla\cap D} .
\end{equation}
Let $\Gamma_\mu$ and $\Gamma_\nu$ be the $\Gamma$-potentials of $\mu$ and $\nu$ respectively, i.e., 

\begin{equation}\label{15} \Gamma_\mu(z)= \int_\erreu \Gamma(\zeta,z) \ d\mu(\zeta)= \unor\int_\omegapalla \Gamma (\zeta,z) W(z_0\menou\circ\zeta)\ d\zeta,
\end{equation} and 
\begin{equation}\label{16} \Gamma_\nu(z)= \int_\erreu \Gamma(\zeta,z) \ d\nu(\zeta)= \unor\int_D \Gamma (\zeta,z) W(z_0\menou\circ\zeta)\ d\zeta.
\end{equation} 
From \eqref{otto} and assumption \eqref{trestelle}, we have 

\begin{equation}\label{17} \Gamma_\mu(z)= \Gamma(z_0,z)\quad \forall\, z\notin \omegapalla
\end{equation} and 

\begin{equation}\label{18} \Gamma_\nu(z)= \Gamma(z_0,z)\quad \forall\, z\notin D.
\end{equation} 
Therefore, 
\begin{equation}\label{19} \Gamma_\mu(z)= \Gamma_\nu(z)\quad \forall\, z\in \erreu \meno (\omegapalla\cup D).
\end{equation}  
Let us remark that 
\begin{equation*}(\Gamma_\mu-\Gamma_\nu)(z) =\int_\erreu \Gamma(z,\zeta) f(\zeta) \ d\zeta, \end{equation*} 
where
$$ f(\zeta):=\unor (\mathbb{1}_\omegapalla -\mathbb{1}_D)(\zeta) W(z_0\menou\circ\zeta).$$
From assumption \eqref{nove}, the function $f\in L^p(\erreu)$ for a suitable $p>\frac{Q}{2}$. Furthermore, as $\omegapalla$ and $D$ are bounded, $f$ has compact support. Then, 

\begin{equation}\label{20} \Gamma_\mu - \Gamma_\nu \in C(\erreu,\erre)
\end{equation} 
(see Subsection \ref{subsection 2.7}). Then, since $\Gamma\mu$ is finite at any point (see \eqref{stellina}), we have 
$$\Gamma_\nu(z)<\infty \quad\forall\, z\in\erreu.$$

%
%
%

With all these results at hand, we can prove Theorem \ref{main} with a procedure inspired by the one used in \cite{cuplanidea}. We split our procedure into several steps in which we simply denote by $\Omega$ the $\elle$-ball $\omegapalla.$ 

\stepone\label{stepone} The aim of this step is to prove the inequality

\begin{equation}\label{22} \Gamma_\mu \le \Gamma_\nu \inn \erreu \setminus \Omega.\end{equation} 

To this end, keeping in mind \eqref{19} and remarking that 

$$\erreu\meno\Omega \subseteq (\erreu\meno(\Omega \cup D))\cup D,$$

it is enough to prove that 

\begin{equation}\label{23} \Gamma_\mu (z) \le \Gamma_\nu (z) \quad \forall\, z\in D.\end{equation} 

Define 

\begin{equation}\label{24} \mu_1:= \mu|_{\Omega \meno D} \andd \nu_1:= \nu|_{D \meno \Omega}.\end{equation} 
Then $\mu_1$ and $\nu_1$ are compactly supported Radon measures whose $\Gamma$-potentials are finite at any point of $\erreu$.
If $z\in D$ we have 
\begin{eqnarray*} (\Gamma_\mu-\Gamma_\nu)(z) &=& \int_{\Omega\cap D} \Gamma(z,\zeta)\ d\mu(\zeta) + \int_{\Omega\meno D} \Gamma(z,\zeta)\ d\mu (\zeta) \\&& - \left(\int_{ D\cap \Omega} \Gamma(z,\zeta)\ d\nu(\zeta) + \int_{D\meno \Omega} \Gamma(z,\zeta)\ d\nu (\zeta)\right),\end{eqnarray*}

so that, keeping in mind that $$\mu|_{\Omega \cap D}=\nu|_{\Omega \cap D}$$

(see \eqref{14}), we get 
\begin{equation}\label{25} (\Gamma_\mu -  \Gamma_\nu ) (z) =  (\Gamma_{\mu_1} -  \Gamma_{\nu_1}) (z) 
 \quad \forall\, z\in D.\end{equation} 
We let 
\begin{equation}\label{26} \hat{u}: = \Gamma_{\mu_1} -  \Gamma_{\nu_1}.\end{equation} 
Since $D$ is open, one has the following inclusions: 
$$D\cap \supp \mu_1 \subseteq D\cap (\overline{\Omega\meno D}) \subseteq D\cap  (\erren \meno D) =\emptyset.$$

As a consequence, 
$$\Gamma_{\mu_1} \mbox{ is $\ellestar$-harmonic in $D$}$$

(see Subsection \ref{subsection 2.5}). It follows that 

$$ \hat{u}\in \underline\elle^\ast (D).$$
We claim that 
\begin{equation}\label{27} \limsup_{D\ni z\rightarrow \zeta} \hat{u}(z) \le 0\quad \forall\, \zeta\in\partial D.\end{equation} 

Taking this claim for granted, for a moment, from the Maximum Principle for $\ellestar$-subharmonic functions (see Subsection \ref{subsection 2.3}), we get 
$$\hat{u}\le 0\inn D.$$
Then, by \eqref{25} and \eqref{26}, 
$$\Gamma_\mu-\Gamma_\nu \le 0 \inn D.$$ 
This proves \eqref{23}, hence \eqref{22}.

We are left with the proof of claim   \eqref{27}. Keeping again in mind  \eqref{26}  and \eqref{25}, and using the continuity of 
the function $\Gamma_\mu-\Gamma_\nu$ - see  \eqref{20} - for every $\zeta \in \partial D$ we have

\begin{eqnarray}\label{28} \limsup_{D\ni z\rightarrow \zeta} \hat{u}(z)  &= &  \limsup_{D \ni z\rightarrow \zeta} (\Gamma_{\mu_1}(z) - \Gamma_{\nu_1}(z))= \limsup_{D \ni z\rightarrow \zeta} (\Gamma_\mu(z) - \Gamma_\nu(z)) \\ \nonumber &=&  \lim_{z\rightarrow \zeta} (\Gamma_\mu(z) - \Gamma_\nu(z)) =
 \Gamma_\mu(\zeta) - \Gamma_\nu(\zeta). \end{eqnarray} 
On the other hand,  by \eqref{stellina},  $ \Gamma_\mu(\zeta) \le \Gamma(z_0,\zeta)$, while, since $\zeta \in \partial D$ so that $\zeta\notin D$, from \eqref{18} it follows $\Gamma_\nu(\zeta)=\Gamma(z_0,\zeta).$ Using this information in  \eqref{28} we finally get 
 \begin{eqnarray*} \limsup_{D\ni z\rightarrow \zeta} \hat{u}(z)  \le 0. \end{eqnarray*} 
 
 This proves the claim \eqref{27} and completes the proof of \eqref{22}.

 \steptwo\label{steptwo}
 The aim of this step is the proof of the inclusion 
 \begin{equation}\label{29} D\subseteq \Omega.\end{equation} 
Since $\supp \mu =\overline \Omega$, the potential function $\Gamma_\mu$ is $\ellestar$-harmonic in $\erreu\meno\overline\Omega.$ Therefore
 
  \begin{equation}\label{30} v:=\Gamma_\mu  -\Gamma_\nu \in \underline\elle^\ast(\erreu\meno \overline\Omega) .\end{equation} 
We know, by {\it \hyperref[stepone]{Step I}},  that
 \begin{equation}\label{31} v\le 0 \inn \erreu\meno \overline\Omega .\end{equation} 
 On the other hand, by \eqref{19},
 
  \begin{equation}\label{32} v\equiv 0 \inn \erreu\meno (\Omega \cup D) .\end{equation} 
  
  Now, let $z=(x,t)$ be an arbitrary point of $\erreu\meno \overline\Omega$. From the very definition of $\Omega$ ($=\omegapalla$ and $z_0=0$) the bounded  subset of $\erren$ 
    \begin{equation*} \Omega_t:= \left\{ \xi\in\erren\ |\ (\xi,t)\in\overline\Omega\right\}
    \end{equation*}
  is empty or convex. As a consequence, for every fixed $j\in \{1, \ldots, p_0\},$
  
   \begin{equation*} z + s e_j \in \erreu\meno\overline\Omega\quad\forall\, s\geq 0
    \end{equation*} or
     \begin{equation*} z - s e_j \in \erreu\meno\overline\Omega\quad\forall\, s\geq 0.
    \end{equation*}
  To fix ideas, let us suppose that the first case occurs. Then, since $\Omega\cup D$ is bounded, there exists $S>0$ such that 
  \begin{equation*} z^\ast:=z +  S e_j \in \erreu\meno (\Omega\cup D) 
    \end{equation*} and 
   \begin{equation*} z + s e_j \in \erreu\meno \overline\Omega \quad \forall\, s\in [0,S].
    \end{equation*} 
    Then, by \eqref{32} and \eqref{31}, 
      \begin{equation*} v(z^\ast)=0= \max_{\erreu \meno \overline \Omega} v.
    \end{equation*} 
    The propagation of maxima for $\ellestar$-subharmonic functions (see Subsection \ref{subsection 2.4}) implies that

  \begin{equation*} v(z + s e_j)=v(z^\ast)=0  \quad \forall\, s\in [0,S].
    \end{equation*} 
    In particular, for $s=0$, we get $v(z)=0.$ Since $z$ is an arbitrary point of $\erreu\meno\overline\Omega,$ we have so proved that 
     \begin{equation*} v\equiv 0 \inn \erreu\meno\overline\Omega.
    \end{equation*} 
    This means that
     \begin{equation*} \Gamma_\mu=\Gamma_\nu  \inn \erreu\meno\overline\Omega.
    \end{equation*} 
On the other hand, as we have already observed, $\Gamma_\mu$ is $\ellestar$-harmonic in  $\erreu\meno\overline\Omega.$
As a consequence, 
 \begin{equation*} \nu= - \ellestar (\Gamma_\nu)=- \ellestar (\Gamma_\mu) = 0   \inn \erreu\meno\overline\Omega, 
    \end{equation*} i.e., 
    
   \begin{equation*} \nu (\erreu\meno\overline\Omega)=0, 
    \end{equation*} or, equivalently, 
       \begin{equation*} \supp \nu \subseteq \overline \Omega.
    \end{equation*} 
    Then $D\subseteq \overline \Omega$ since $\supp \nu=\overline D.$ 
    As a consequence, $D\subseteq \mathrm{int }(\overline\Omega).$ On the other hand, keeping in mind that $\gamma$ is $\delta_\lambda$-homogeneous of degree $2-Q$, it is easy to show that $$\mathrm{int }(\overline \Omega)=\Omega.$$
   Hence, from the last inclusion, we get
  $$D\subseteq \Omega.$$
    \stepthree In this final step we prove that
    \begin{equation}\label{33} D=\Omega .\end{equation} 
    We argue by contradiction and assume $D\neq\Omega.$ In this case, since $D\subseteq\Omega$ 
    by {\it \hyperref[steptwo]{Step II}}, 
there exists $z\in\Omega$ such that $z\notin D$. Then, by inequality \eqref{undici}
    
      \begin{equation}\label{34} \Gamma(z_0,z)>\Gamma_\mu(z) .\end{equation} 
On the other hand, since $D\subseteq \Omega$, from \eqref{14} we get 
 \begin{equation}\label{35} \mu|_D=\nu .\end{equation} 
 Moreover, since $z\notin D$, 
 \begin{equation}\label{strano} \Gamma_\nu (z)=\Gamma (z_o,z).\end{equation} 
Putting \eqref{34}, \eqref{35} and \eqref{strano} together we have 
    
    \begin{eqnarray*} \Gamma(z_0,z) &>&\Gamma_\mu(z)=\int_\Omega \Gamma(\zeta,z)\ d\mu(\zeta)\\ 
&\geq& \int_D \Gamma(\zeta,z)\ d\mu(\zeta)  = \int_D \Gamma(\zeta,z)\ d\nu(\zeta) \\ &=&\Gamma_\nu(z)=\Gamma(z_o,z),     \end{eqnarray*}
    that is $\Gamma(z_0,z)>\Gamma(z_0,z).$ This contradiction proves \eqref{33} and completes the proof of Theorem \ref{main}.

\section{Proof of Corollary \ref{cor}}

We first remark that Corollary's assumptions imply 
    \begin{equation}\label{37} D\subseteq \erren\times ]-\infty,t_0[,\end{equation} 
    if $t_0$ is the time component of $z_0$ (i.e. $z_0=(x_0,t_0)$ for a suitable $x_0\in\erren$). Indeed, assume by contradiction that
    
     \begin{equation*}\label{36} D\cap (\erren\times [t_0,\infty[ ) \neq \emptyset.\end{equation*} 
    Then, since $D$ is bounded, there exists $z^\ast=(x^\ast,t^\ast)\notin D$ such that $,t^\ast>t_0$ and 
    
    $$\dstar:=D\cap ( \erren\times ]\tstar,\infty[) \neq \emptyset.$$
Let us now consider the function 
$$\ustar(\zeta):=\Gamma(\zeta,\zstar),\ \zeta\in\erreu.$$
Since $\zstar \notin D, \ustar$ is $\elle$-harmonic and nonnegative in $D$. Moreover 

$$\ustar(z_0)=\Gamma(z_0,\zstar)=0\andd \ustar>0\inn \dstar.$$
Then, by assumption \eqref{star!} 

\begin{equation*} \ustar(z_0)= \unor \int_{D}  \ustar(\zeta) W(z_0^{-1} \circ \zeta)\ d \zeta .\end{equation*}
This is a contradiction since $\ustar(z_0)=0$ and 

\begin{equation*}  \int_{D}  \ustar(\zeta) W(z_0^{-1} \circ \zeta)\ d \zeta \geq   \int_{\dstar}  \ustar(\zeta) W(z_0^{-1} \circ \zeta)\ d \zeta >0. \end{equation*}
This contradiction proves the inclusion \eqref{37}. 

Let us now observe that, for every $z\notin D, z\neq z_0,$ the function $$\zeta\frecciaf \Gamma(\zeta,z)$$
is non-negative and $\elle$-harmonic in an open set, precisely $\erreu\meno\{z\}$, containing $D\cup\{z_0\}.$ Then by assumption \eqref{star!} 

\begin{equation} \label{38} \int_{D}  \Gamma(\zeta,z) W(z_0^{-1} \circ \zeta)\ d \zeta = r\Gamma (z_0,z)\end{equation}
for every $z\notin D,\ z\neq z_0.$ On the other hand, identity \eqref{38} is trivial if $z=z_0$ since $\Gamma(z_0,z_0)=0$ and, being $D\subseteq \erren\times ]-\infty, t_0[,$ $\Gamma(\zeta,z_0)=0$ for every $\zeta\in D$ so that 
\begin{equation*}  \int_{D}  \Gamma(\zeta,z) W(z_0^{-1} \circ \zeta)\ d \zeta = 0.\end{equation*}
We have so proved that $D$ satisfies the assumption of Theorem \ref{main}, keeping in mind that we are assuming condition \eqref{nove}. Then, by Theorem \ref{main},
$$D=\omegapalla.$$

\subsection*{Acknowledgment}
The first author  has been partially supported by the Gruppo Nazionale per l'Analisi Matematica, la Probabilit\`a e le
loro Applicazioni (GNAMPA) of the Istituto Nazionale di Alta Matematica (INdAM).
\bibliographystyle{plain}

\end{document}